\documentclass[12pt]{article}
\usepackage{amsmath}
\usepackage{amssymb}
\newsymbol \blackbox 1004
\newcommand{\eh}{\hfill}\newlength{\sperr}

\newenvironment{proof}{{\settowidth{\sperr}{\bf\rm
Proof}%
\par\addvspace{0.3cm}\noindent\parbox[t]{1.3\sperr}
{\bf\rm P\eh r\eh o\eh o\eh f\eh }%
}}{\nopagebreak\mbox{}
$\blackbox$\par\addvspace{0.3cm}}

\def\a{\alpha}
\def\b{\beta}

\def\ov{\overline}

\def\vp{\varphi}

\def\wh{\widehat}
\def\wt{\widetilde}
\def\ov{\overline}
\def\p{\partial}
\def\BC{{\mathbb C}}

\newtheorem{Pa}{Paper}[section]
\newtheorem{Tm}[Pa]{{\bf Theorem}}

\newtheorem{Cy}[Pa]{{\bf Corollary}}
\newtheorem{Rk}[Pa]{{\bf Remark}}

\newtheorem{Ee}[Pa]{{\bf Example}}
\newtheorem{Dn}[Pa]{{\bf Definition}}
\newtheorem{Pn}[Pa]{{\bf Proposition}}
\newcommand{\CC}
{{\mathchoice {\setbox0=\hbox{$\displaystyle\rm
C$}\hbox{\hbox
to0pt{\kern0.4\wd0\vrule height0.9\ht0\hss}\box0}}
{\setbox0=\hbox{$\textstyle\rm C$}\hbox{\hbox
to0pt{\kern0.4\wd0\vrule height0.9\ht0\hss}\box0}}
{\setbox0=\hbox{$\scriptstyle\rm C$}\hbox{\hbox
to0pt{\kern0.4\wd0\vrule height0.9\ht0\hss}\box0}}
{\setbox0=\hbox{$\scriptscriptstyle\rm C$}\hbox{\hbox
to0pt{\kern0.4\wd0\vrule height0.9\ht0\hss}\box0}}}}

\title{Skew-self-adjoint discrete and continuous Dirac
type systems:
inverse problems and Borg-Marchenko theorems}

\author{Alexander Sakhnovich}

\date{}
\begin{document}
\maketitle

{\bf Short title.}  Skew-self-adjoint Dirac type
systems
\vspace{5mm}

Murrhardter Str. 46, 73642 Welzheim, Germany, \\ e-mail address:
al$_-$sakhnov@yahoo.com

\begin{abstract}
New formulas on the inverse problem for the continuous
skew-self-adjoint Dirac type system are obtained. For the discrete
skew-self-adjoint Dirac type system the solution of a general type
inverse spectral problem  is also derived in terms of the Weyl
functions. The description of the Weyl functions on the interval
is given. Borg-Marchenko type uniqueness theorems are derived for
both discrete and continuous non-self-adjoint systems too.

\end{abstract}

\section{Introduction} \label{intro}
\setcounter{equation}{0}
Skew-self-adjoint Dirac type system
\begin{equation}       \label{0.1}
\frac{d}{dx}u(x, \lambda )=\Big(i \lambda j+j
V(x)\Big)u(x,
\lambda ), \quad x \geq 0,
\end{equation}
where
\begin{equation}   \label{0.2}
j = \left[
\begin{array}{cc}
I_{p} & 0 \\ 0 & -I_{p}
\end{array}
\right], \hspace{1em} V= \left[\begin{array}{cc}
0&v\\v^{*}&0\end{array}\right],
 \end{equation}
$I_p$ is the $p \times p$ identity matrix and $v$ is
$p \times p$
matrix function, is a classical object of analysis. It
is also
called Zakharov-Shabat or AKNS. In particular, system
(\ref{0.1})
and its discrete analog
\begin{equation} \label{0.3}
W_{k+1}( \lambda)-W_k(\lambda)=- \frac{i}{\lambda}C_k W_k(\lambda), \quad
C_k=C_k^*=C_k^{-1},\quad k=0, \, 1, \ldots
\end{equation}
are auxiliary linear systems for many important integrable
non-linear equations \cite{AS, FT, KS}. Various results and
references on the scattering theory for system (\ref{0.1}) one can
find in \cite{FT}. Weyl functions for this system on the interval
and semi-axis have been introduced in \cite{SaA1}. (Weyl functions
are also called Weyl-Titchmarsh or $M$-functions in the
literature.) The existence and uniqueness of the Weyl function for
system (\ref{0.1}) with a bounded on the semi-axis potential was
proved in \cite{SaA1}, and inverse problems in terms of the Weyl
functions on the interval and semi-axis have been solved.
Interesting recent spectral results on the non-self-adjoint
(especially skew-self-adjoint) Dirac type systems one can find in
\cite{GL, GrK, KlS}. See also further developments of the
Weyl-Titchmarsh theory with respect to system (\ref{0.1}) in
\cite{CG2, GKS0, SaA2}. Weyl functions are successfully used in
solving initial-boundary value problems for integrable nonlinear
equations (see \cite{B, Kv, SaA2, SaA4, SaL1, SaL2, SaL3}).
Moreover, Borg-Marchenko type results on the potentials coinciding
on a part of the interval, where the systems are defined, in terms
of the asymptotics of the Weyl functions of these systems are of
great current interest, and the self-adjoint case have been
considered recently in \cite{CG, CGR, GS, GZ, SaA3, Si, Si2}.

In Section 2 we shall obtain Borg-Marchenko type result and some
new formulas on the inverse problem  for skew-self-adjoint system
(\ref{0.1}).

The largest Section 3 of the paper is dedicated to the discrete
case - system (\ref{0.3}). In the case of the discrete system
(\ref{0.3}) on the semi-axis explicit procedure to recover
potential  from the rational Weyl function one can find in
\cite{KS}. Here we shall solve inverse problem for the discrete
case on the interval without requiring the Weyl function to be
rational. The set of the Weyl functions will be completely
described in terms of the Taylor coefficients. Borg-Marchenko type
theorem will be derived too.

We denote the complex plane by $\BC$, upper (lower) semi-plane by
$\BC_+$ ($\BC_-$),  and we denote the set of bounded operators
acting from ${\bf H}_1$ into ${\bf H}_2$ by $\{ {\bf H}_1, \, {\bf
H}_2 \}$.

%%%%%%%%%%%%%%%%%%%%%%%%%%%%%%%%%%%%%%%%%%%%%%%%%%%%%%%%%%%%%%%%%%%%
\section{Continuous case} \label{Cont}
\setcounter{equation}{0}

In this paper we shall consider systems (\ref{0.1}) with locally
bounded potentials on the intervals $[0, \, l]$ and $[0, \,
\infty)$.
 In other words we assume that the inequalities
\begin{equation} \label{0.4}
\| v(x) \| \leq M \quad (0 \leq x \leq l)
\end{equation}
are true. Normalize  the $m \times m$ ($m=2p$)
fundamental
solution $u$ of system (\ref{0.1})  by the condition
$u(0,\lambda)=I_m$. In view of (\ref{0.1}) it is immediate
that
\begin{equation} \label{0.4.1}
u(x, \lambda)^*u(x, \lambda)=u(x, \lambda)u(x, \lambda)^*=I_m \quad
{\mathrm{for}}
\quad \lambda = \ov{\lambda}.
\end{equation}
We shall use also the notations for the rows of
$u(x,0)$:
\begin{equation} \label{0.4.2}
\b(x)=[I_p \quad 0]u(x, 0), \quad \chi(x)=[0 \quad
I_p]u(x, 0).
\end{equation}
From (\ref{0.1}), (\ref{0.4.1}), and (\ref{0.4.2}) one
easily gets
\begin{equation} \label{0.4.3}
\b(0)=[I_p \quad 0], \quad \b \b^* \equiv I_p, \quad
\b^{\prime}
\b^* \equiv 0 \quad (\b^{\prime}=\frac{d}{dx} \b),
\end{equation}
\begin{equation} \label{0.4.4}
\b \chi^* \equiv 0, \quad \chi(0)=[0 \quad I_p], \quad
\chi \chi^*
\equiv I_p, \quad \chi^{\prime} \chi^* \equiv 0,
\end{equation}
\begin{equation} \label{0.4.5}
v(x)=\b^{\prime}(x) \chi(x)^*.
\end{equation}
Put now
\begin{equation} \label{0.5}
{\cal W}(\lambda)= \{{\cal W}_{ij}(\lambda) \}_{i,j=1}^2=u(l, \ov{\lambda})^*,
\end{equation}
where ${\cal W}_{ij}$ are $p \times p$ blocks of ${\cal W}$.
According to (\ref{0.1}), (\ref{0.4}), and (\ref{0.5}) we have
\cite{SaA1}:
\begin{equation} \label{0.5'}
{\cal W}(\lambda)j {\cal W}(\lambda)^* \leq j, \quad {\cal
W}(\lambda)^*j {\cal W}(\lambda) \leq j \quad {\mathrm{for}} \quad
\Im \lambda \leq -M.
\end{equation}
A pair of $p \times p$ meromorphic in the semi-plane
$\Im \lambda < -M$
matrix functions $R$ and $Q$ is called non-singular
with
$j$-property if
\begin{equation}\label{0.6}
R(\lambda )^*R(\lambda )+Q(\lambda )^*Q(\lambda )>0, \quad
\left[\begin{array}{lr}
  R(\lambda )^* & Q(\lambda )^* \end{array}\right]\, j \,
\left[\begin{array}{c} R(\lambda ) \\  Q(\lambda )
\end{array}\right]\leq 0.
\end{equation}
\begin{Dn} \label{Dn1.1} \cite{SaA1}
Let system (\ref{0.1}) be given on the interval $[0,
\, l]$ and
satisfy (\ref{0.4}). Then the linear-fractional
transformations
\begin{equation}\label{0.7}
\varphi (\lambda ) = \bigl( {\cal W}_{11}(\lambda )R(\lambda )+ {\cal
W}_{12}(\lambda
)Q(\lambda )\bigr) \bigl( {\cal W}_{21}(\lambda ) R(\lambda )+{\cal
W}_{22}(\lambda
)Q(\lambda )\bigr)^{-1},
\end{equation}
where $\Im \lambda < -M$ and $R$, $Q$  are non-singular
pairs with
$j$-property, are called Weyl functions of this
system.
\end{Dn}
We denote the set of Weyl functions by ${\mathcal N}(l)$. From
(\ref{0.5'})-(\ref{0.7}) it follows that
\begin{equation}\label{0.8}
\vp(\lambda)^* \vp(\lambda) \leq I_p.
\end{equation}
The procedure to recover a system from its Weyl functions can be
given in terms of the $S$-nodes introduced in
\cite{SaL0}-\cite{SaL3} (see also references therein). Namely,
operator $S=S(l)$ acting in the space $L^2_p(0,l)$  of squarely
integrable vector functions and block operator $\Pi= [\Phi_1 \quad
\Phi_2]$, where $\Phi_k \,$ ($k=1,2$) act from $\BC^p$ into
$L^2_p(0,l)$, are defined in \cite{SaA1} and satisfy operator
identity
\begin{equation}\label{0.8'}
AS-SA^*=i \Pi \Pi^*, \quad A \in \{L^2_p(0,l), \, L^2_p(0,l) \},
\quad A=i \int_0^x \, \cdot \, dt.
\end{equation}
The construction of these operators is based on the relations:
\begin{equation}\label{2n}
S=V_-^{-1}\Big(V_-^{-1} \Big)^*, \quad \Phi_k=V_-^{-1} \b_k, \quad
K:=i \b(x) \int_0^x \b(t)^* \, \cdot \, dt=V_-AV_-^{-1},
\end{equation}
where $V_-$ is a triangular and bounded in $L^2_p(0, \, l)$
together with its inverse operator, $\b_k$ are $p \times p$ blocks
of $\b$. (Here "$V_-$ is triangular" means that
$P_xV_-=P_xV_-P_x$, $0<x<l$, where $P_x$ is the orthogonal
projector from $L^2_p(0,l)$ onto $L^2_p(0,x)$.) Moreover operator
$V_-^{-1}$ takes functions with bounded derivatives into functions
with bounded derivatives, $\Big(V_-^{-1}f \Big)(0)=f(0)$, and
$V_-$ is normalized so that $\Phi_1$ proves a natural embedding:
$\Phi_1 g \equiv g = {\mathrm{const}}$. Operator $\Phi_2$ given in
the second relation in (\ref{2n}) can be presented as a
multiplication by matrix function and we denote this matrix
function by $s$: $\, \Phi_2 g =s(x) g$. It proves \cite{SaA1} that
the following asymptotics takes place for the Weyl functions of
system (\ref{0.1}) on $[0, \, l]$:
\begin{equation}\label{2n'}
{\displaystyle \vp(\lambda)=-2i \lambda \int_0^le^{-2i \lambda x }s(x)dx+O(|\lambda
e^{-2i \lambda l }|). }
\end{equation}
Using (\ref{2n'})
 a $p \times p$ matrix function $s(x)$ with
the entries from $L^2(0,l)$ is recovered in \cite{SaA1} via the
Fourier transform
\begin{equation}\label{0.9}
e^{- \eta x}s(x)= \frac{i}{2 \pi} \int_{- \infty}^{\infty}e^{i \xi
x} \lambda^{-1} \vp(\lambda /2) d \xi \quad (\lambda= \xi - i \eta , \quad \eta
>2M),
\end{equation}
where $s$ does not depend on the choice of $\eta >2M$.
\begin{Tm} \label{Tm1.2} \cite{SaA1} Suppose $\vp$ is
a Weyl function of system (\ref{0.1}) on $[0, \, l]$. Then system
(\ref{0.1}) is uniquely recovered from $\vp$ by the following
procedure. First use (\ref{0.9}) to introduce operators $\Phi_k \,
\in \{ \BC^p, \, L^2_p(0,l) \} $ $(k=1,2):$
\begin{equation}\label{0.10}
\Phi_1 g \equiv g ={\mathrm{const}}, \quad \Phi_2 g =s(x)g,
\end{equation}
and operator $S=S(l)$ acting in $L^2_p(0,l)$:
\begin{equation}\label{0.11}
 S f=\frac{1}{2} \frac{d}{d x}\int_0^l \left(
\frac{\p}{\p t} \int^{x+t}_{|x-t|} \left( I_p+s \Big(
\frac{r+x-t}{2} \Big) s \Big( \frac{r+t-x}{2} \Big)^* \right)d r
\right)f(t) d t.
\end{equation}
These operators are well-defined and bounded.
Moreover, $S$ is
positive, i.e., $S>0$, and $S$ is boundedly
invertible.

Next we recover $\b^* \b$ by the formula
\begin{equation}\label{0.11'}
\b(x)^*\b(x)= \frac{d}{d x} \Big( \Pi^* S(x)^{-1} P_x \Pi \Big),
\end{equation}
where $\Pi$ is a block operator $\Pi:= [\Phi_1 \quad \Phi_2]$, $\,
P_x$ is the orthogonal projector from $L^2_p(0,l)$ onto
$L^2_p(0,x)$, and operator $S(x) \in \{L^2_p(0,x), \, L^2_p(0,x)
\}$ is given by the equality $S(x)=P_x S(l) P_x$.

Finally, the potential $v$ is uniquely recovered from
$\b^* \b$
using relations (\ref{0.4.3})-(\ref{0.4.5}).
\end{Tm}
Notice that $P_x \Pi $ in (\ref{0.11'}) is considered as a matrix
function $[I_p \quad s(t)]$ and $S(x)^{-1}$ is applied to this
matrix function columnwise.

There is a simpler way to recover $\chi$ from $\vp$ so that $\b$
and $v$ could be recovered after that from $\chi$.
\begin{Pn} \label{Pn3} Suppose $\vp$ is a Weyl function of system (\ref{0.1})
with a bounded on $[0, \, l]$ potential $v$. Then matrix function
$\chi$ defined by the second relation in (\ref{0.4.2}) is given by
the equality
\begin{equation}\label{0.12}
\chi(x)=[0 \quad I_p]- \int_0^x \Big(S(x)^{-1} s'(t)\Big)^*[I_p
\quad s(t)]dt.
\end{equation}
\end{Pn}
\begin{proof}. Denote the upper bound of $||v||$ by $M$, i.e.,
assume that (\ref{0.4}) holds. As $\Big(V_-^{-1}f \Big)(0)=f(0)$
according to the first relation in (\ref{0.4.2}) and second
relation in (\ref{2n}) we have:
\begin{equation}\label{0.13}
s(+0)=\Big(V_-^{-1} \b_2 \Big)(0)= \b_2(0)=0.
\end{equation}
Recall that $\b$ has a bounded derivative and so $s=V_-^{-1} \b_2$
has a bounded derivative too. Apply now operators on the both
sides of the third relation in (\ref{2n}) columnwise to the matrix
function $V_-s'(x)$. One easily gets
\begin{equation}\label{3n}
V_-s(x)=\b(x) \int_0^x \b(t)^* \Big(V_-s'\Big)(t)dt.
\end{equation}
From the second relation in (\ref{2n}) and from (\ref{3n}) it
follows that
\begin{equation}\label{4n}
{\displaystyle \b(x) \left( \left[
\begin{array}{c}
0 \\ I_{p}
\end{array}
\right] \, - \, \int_0^x \left[
\begin{array}{c}
I_{p} \\ s(t)^*
\end{array}
\right]\left(V_-^*P_xV_-s'\right)(t)dt\right)=0.}
\end{equation}
As according to the first relation in (\ref{2n}) we have
$S(x)^{-1}=V_-^*P_xV_-$ formula (\ref{4n}) implies $\b(x) \wh
\chi(x)^*=0$, where $ \wh \chi$ denotes expression on the
right-hand side of (\ref{0.12}). Using $S(x)^{-1}=V_-^*P_xV_-$ we
can rewrite the right-hand side of (\ref{0.12}) as
\[
\wh \chi(x)=[0 \quad I_p]- \int_0^x \Big(V_- s'\Big)(t)^* \b(t)dt.
\]
Thus we have $\wh \chi'(x)=- \Big(V_- s'\Big)(x)^* \b(x)$, which
taking into account $\b \wh \chi^*=0$ implies $\wh \chi' \wh
\chi^*=0$. It is immediate also that $\wh \chi(0)=[0 \quad I_p]$.
Compare equalities
\[
\b(x) \wh \chi(x)^*=0, \quad \wh \chi(0)=[0 \quad I_p], \quad \wh
\chi'(x) \wh \chi(x)^*=0
\]
with the equalities (\ref{0.4.4}) to see that $\chi \equiv \wh
\chi$, i.e., (\ref{0.12}) is true.
\end{proof}
For the case $p=1$ formula (\ref{0.12})  have been announced in
\cite{SaA2}. In that case $\b$ is recovered from $\chi=[\chi_1
\quad \chi_2]$ in the easiest way: $\b=[\ov \chi_2 \quad - \ov
\chi_1]$.

Notice also that as matrix function $s(x)$ has a bounded
derivative and $s(0)=0$, formula (\ref{0.11}) can be written down
in a simpler way:
\begin{equation}\label{0.14}
 S f=f(x)+ \frac{1}{2}\int_0^l
 \int^{x+t}_{|x-t|} s' \Big(
\frac{r+x-t}{2} \Big) s' \Big( \frac{r+t-x}{2} \Big)^*d r f(t) d
t.
\end{equation}

Finally, let us formulate Borg-Marchenko type theorem.
\begin{Tm} \label{Tm4}
Let $\varphi_1$ and $ \varphi_2$ be the Weyl functions of the two
Dirac type systems (\ref{0.1}) on $[0, \, l]$ with the
 bounded potentials denoted by $v_1$ and $v_2$,
respectively. Suppose that on some ray $\Im \lambda=c \Re \lambda <0$ and
for some $0<r<2l$   we have
\begin{equation} \label{0.15}
||\varphi_1 (\lambda)-  \varphi_2(\lambda)||= O(\lambda e^{-i \lambda r}) \quad (| \lambda
| \to \infty).
\end{equation}
Then $v_1(x) \equiv v_2(x)$ for $0<x< \frac{r}{2}$.
\end{Tm}
\begin{proof}. From (\ref{2n'}) it follows that for
each $r<2l$, $| \lambda | \to \infty$, and $ \Im \lambda / | \Re \lambda | < -
\delta < 0 $ we have
\begin{equation} \label{0.16}
|| \varphi_1 (\lambda)-  \varphi_2(\lambda)||= - i \lambda \int_0^re^{-i \lambda
t}\Big( s_1 \big(\frac{t}{2}\big)-
s_2\big(\frac{t}{2}\big)\Big)dt+o(\lambda e^{-i \lambda r}),
\end{equation}
 where $s=s_1$ and $s=s_2$ correspond to systems
(\ref{0.1}) with $v=v_1$ and $v=v_2$, respectively. Formulas
(\ref{0.15}) and (\ref{0.16}) imply the equality
\begin{equation} \label{0.17}
F(\lambda):= \int_0^re^{-i \lambda (t-r)}\Big( s_1 \big(\frac{t}{2}\big)-
s_2\big(\frac{t}{2}\big)\Big)dt=O(1)
\end{equation}
on the ray $\Im \lambda=c \Re \lambda$. Notice that $F(\lambda)$ is bounded in
the closed upper semi-plane $\Im \lambda \geq 0$.  Using now
Phragmen-Lindel\"of theorem for an angle, in view of  (\ref{0.17})
we derive that $F(\lambda)$ is bounded also in $\BC_-$ and thus in the
whole plane. Moreover, $F(\lambda) \to 0$ on the rays in $\BC_+$, i.e.,
$F(\lambda) \equiv 0$. It is immediate that  $s_1(x) \equiv s_2 (x)$
$(x<r/2)$.  Applying the procedure to solve the inverse problem as
in Theorem \ref{Tm1.2} (one can use also formulas (\ref{0.12}),
(\ref{0.14})) we finally get $v_1(x) \equiv v_2 (x)$ for $x<r/2$.
\end{proof}

Consider now system (\ref{0.1}) on the semi-axis $[0, \, \infty)$.
\begin{Dn} \label{Dn5}
Let system (\ref{0.1}) be given on the semi-axis $[0, \, \infty)$.
Then $p \times p$ matrix function $\vp(\lambda)$ analytic  in some
semi-plane $\Im \lambda < -M$ is called Weyl function of this system if
inequalities
\begin{equation} \label{0.18}
\int_0^\infty \left[ \begin{array}{lr}    \varphi (\lambda)^* & I_p
\end{array} \right]
   u(x, \lambda)^*
 u(x, \lambda)
 \left[ \begin{array}{c}
  \varphi (\lambda) \\ I_p \end{array} \right] dx < \infty
\end{equation}
 hold for all $\lambda$ in the semi-plane.
 \end{Dn}
\begin{Pn} \label{Pn6} \cite{SaA1}
System (\ref{0.1}) with bounded on the semi-axis $[0, \, \infty)$
potential $v$
\begin{equation} \label{0.19}
\| v(x) \| \leq M \quad (0 \leq x < \infty)
\end{equation}
has a unique Weyl function $\vp$. Moreover, the matrix disks
${\mathcal N}(l)$ of Weyl functions on the intervals converge to a
point, i.e., to an unique function and this function proves to be
the Weyl function on the semi-axis $\vp= \bigcap_{l<
\infty}{\mathcal N}(l)$.
\end{Pn}
From Theorem \ref{Tm4} and Proposition \ref{Pn6} follows
\begin{Cy} \label{Cy7}
Let $\varphi_1$ and $ \varphi_2$ be the Weyl functions of the two
Dirac type systems (\ref{0.1}) on $[0, \, \infty)$ with the
bounded potentials denoted by $v_1$ and $v_2$, respectively.
Suppose that on some ray $\Im \lambda=c \Re \lambda <0$ and for some $0<r$
equality (\ref{0.15}) is true. Then $v_1(x) \equiv v_2(x)$ for
$0<x< \frac{r}{2}$.
\end{Cy}
\section{Discrete case} \label{Disc}
\setcounter{equation}{0} In this section we shall consider
skew-self-adjoint matrix discrete Dirac type system on the
interval:
\begin{equation} \label{2.1}
W_{k+1}( \lambda)-W_k(\lambda)=- \frac{i}{\lambda}C_k W_k(\lambda), \quad
C_k=C_k^*=C_k^{-1},\quad 0 \leq k \leq n.
\end{equation}
When $p=1$, then either $C_k=\pm I_2$ or  $C_k=U(k)j U(k)^*$
($U(k) U(k)^*=I_2$).  If $C_k=U(k)j U(k)^*$ ($k \geq 0$), system
(\ref{2.1}) is auxiliary system for isotropic Heisenberg magnet
model. Therefore here we  also assume that $C_k=U(k)j U(k)^*$ and
matrices $U(k)$ are unitary, i.e.,
\begin{equation} \label{2.2}
C_k=I_{2p}-2 \b(k)^*\b(k), \quad \b(k)\b(k)^*=I_{p}, \quad 0 \leq
k \leq n,
\end{equation}
where $\b(k)=[\b_1(k) \quad \b_2(k)]=[0 \quad I_p] U(k)^*$ are $p
\times 2p$ matrices with $p \times p$ blocks $\b_1(k)$, $\b_2(k)$.
Introduce now simple additional conditions:
\begin{equation} \label{2.3}
\det \, \b_1(0) \not= 0, \quad \det \, \b(k-1)\b(k)^* \not= 0,
\quad 0 <k \leq n.
\end{equation}
Similar to the continuous case we shall define Weyl functions  of
the system via M\"obius (linear-fractional) transformation
\begin{equation}\label{2.5}
\varphi (\lambda ) = \bigl( {\cal W}_{11}(\lambda )R(\lambda )+ {\cal W}_{12}(\lambda
)Q(\lambda )\bigr) \bigl( {\cal W}_{21}(\lambda ) R(\lambda )+{\cal W}_{22}(\lambda
)Q(\lambda )\bigr)^{-1},
\end{equation}
where we  put
\begin{equation} \label{2.4}
{\cal W}(\lambda)= \{{\cal W}_{ij}(\lambda) \}_{i,j=1}^2=W_{n+1}(
\ov{\lambda})^*.
\end{equation}
Here $2p \times 2p$ solution $W_k$ of  (\ref{2.1}) is normalized
by the condition $W_0(\lambda)=I_{2p}$ and coefficients ${\cal W}_{ij}$
of the M\"obius  transformation are $p \times p$ blocks of ${\cal
W}$. We shall be interested in the properties of $\vp(\lambda)$ in the
neighborhood of $\lambda=i$. So we require that $R$ and $Q$ are $p
\times p$ matrix functions analytic in the neighborhood of $\lambda=i$
and such that
\begin{equation}\label{2.6}
\det \Big( {\cal W}_{21}(i ) R(i )+{\cal W}_{22}(i )Q(i )\Big)
\not=0.
\end{equation}
Such pairs  $R$, $Q$ always exist as the rows of $[{\cal W}_{21}(i
) \quad {\cal W}_{22}(i )]$ are linearly independent:
\begin{equation}\label{2.7}
{\mathrm {rank}} \, [{\cal W}_{21}(i ) \quad {\cal W}_{22}(i )]=p.
\end{equation}
Indeed, to prove (\ref{2.7}) we can take into account (\ref{2.1})
and $C_k=U(k)j U(k)^*$ and derive
\begin{equation}\label{2.8}
 W(-i )=2^{n+1} \prod_{k=0}^n
 \chi(k)^*\chi(k), \quad \chi(k)=[\chi_1(k) \quad
 \chi_2(k)]=[I_p \quad 0] U(k)^*.
\end{equation}
Notice that by $U(k)^*U(k)=I_{2p}$ we have
\begin{equation}\label{2.8'}
\chi(k) \b(k)^*=0, \quad \chi(k) \chi(k)^*=I_p.
\end{equation}
From (\ref{2.4}) and (\ref{2.8}) it follows that
\[
 [{\cal W}_{21}(i ) \quad {\cal W}_{22}(i )]=
 \]
\begin{equation}\label{2.9}
=2^{n+1}
 \chi_2(0)^*\big(\chi(0)\chi(1)^*\big)\big(\chi(1)\chi(2)^*\big)
 \ldots \big(\chi(n-1)\chi(n)^*\big) \chi(n).
\end{equation}
It remains to show that inequalities (\ref{2.3}) imply:
\begin{equation}\label{2.10}
\det \, \chi_2(0) \not= 0, \quad \det \, \chi(k-1)\chi(k)^* \not=
0, \quad 0 <k \leq n.
\end{equation}
Suppose $\det \, \chi_2(0) = 0$. Then there is a vector $f \not=0$
such that $f^*\chi_2(0) = 0$. So, in view of the first relation in
(\ref{2.8'}) we get $0=f^*\chi(0) \b(0)^*=f^*\chi_1(0) \b_1(0)^*$.
As $\det \, \b_1(0) \not= 0$ it follows that $f^*\chi_1(0)=0$ and
thus $f^*\chi(0)=0$. But according to the second relation in
(\ref{2.8'}) the lines of $\chi(0)$ are linearly independent and
we come to a contradiction. The first inequality in (\ref{2.10})
follows.

Suppose there is a vector $f \not=0$ such that
$f^*\chi(k-1)\chi(k)^* = 0$.  By (\ref{2.8'}) we then obtain
\[
f^*\chi(k-1)= \wt f^* \b(k), \quad \wt f \not=0.
\]
Hence, taking into account the second inequality in (\ref{2.3}) we
have
\[
f^*\chi(k-1) \b(k-1)^*=\wt f^* \b(k) \b(k-1)^* \not= 0.
\]
As the first equality in (\ref{2.8'}) yields $\chi(k-1)
\b(k-1)^*=0$, we get contradiction, i.e., the second inequality in
(\ref{2.10}) is valid too. Now relations (\ref{2.9}) and
(\ref{2.10}) imply equality (\ref{2.7}).

\begin{Dn} \label{Dn2.1}
Let system (\ref{2.1}) be given on the interval $0 \leq k \leq n$
and satisfy (\ref{2.2}), (\ref{2.3}). Suppose $R$ and $Q$ are $p
\times p$ matrix functions analytic in the neighborhood of $\lambda=i$
and such that inequality (\ref{2.6}) holds. Then linear-fractional
transformations $\vp$ of the form (\ref{2.5}) are called Weyl
functions of this system. The pair $R$, $Q$ satisfying our conditions
is called admissible.
\end{Dn}
\begin{Ee} \label{Een1} Put $n=1$, $C_0=-j$, $C_1=J$, where
\begin{equation}\label{e1}
J = \left[
\begin{array}{lr}
0 & I_{p}  \\ I_{p} & 0
\end{array}
\right]=KjK^*, \quad K=\frac{1}{\sqrt 2}\left[
\begin{array}{lr}
I_p & -I_{p}  \\ I_{p} & I_p
\end{array}
\right], \quad K^*=K^{-1}.
\end{equation}
In view of (\ref{e1}) we have $C_0=JjJ$, and so $U(0)=J$. By the second equality in (\ref{e1})
we get $U(1)=K$. It follows that
\begin{equation}\label{e2}
\b(0)=[0 \quad I_p]U(0)^*=[I_p \quad 0], \quad \b(1)=[0 \quad I_p]U(1)^*=\frac{1}{\sqrt 2}[-I_p \quad I_p].
\end{equation}
Thus the conditions (\ref{2.2}) and (\ref{2.3}) are fulfilled. Moreover we have
\[
W_2(\lambda)=\Big(I_{2p} -\frac{i}{ \lambda}J
\Big)\Big(I_{2p} +\frac{i}{ \lambda}j
\Big)=I_{2p}+\frac{i}{ \lambda}(j-J)+\frac{1}{ \lambda^2}Jj, \quad {\mathrm{i.e.,}}
\]
\begin{equation}\label{e3}
{\cal W}(\lambda)=I_{2p}-\frac{i}{ \lambda}(j-J)-\frac{1}{ \lambda^2}Jj.
\end{equation}
In particular, we have ${\cal W}_{21}(i)=2$, ${\cal W}_{22}(i)=2$. Thus condition (\ref{2.6})
takes the form
\begin{equation}\label{e4}
\det \, \big(R(i)+Q(i) \big) \not=0.
\end{equation}
By (\ref{e3}) M\"obius transformation (\ref{2.5}), that describes Weyl functions, takes the form
\begin{equation}\label{e5}
\vp(\lambda)=\frac{\lambda-i}{\lambda+i}(\lambda R+i Q)(i R+ \lambda Q)^{-1}.
\end{equation}
\end{Ee}

Denote the set of Weyl functions of system (\ref{2.1}) 
on the interval $0 \leq k \leq n$
by ${\cal N}(n)$.
\begin{Rk} \label{Rk2.1'} The sets of Weyl functions 
introduced by Definition \ref{Dn2.1} are decreasing. That is, given system (\ref{2.1})
we have ${\cal N}(l) \supseteq {\cal N}(n)$ for $l<n$. Indeed,  the pair
\begin{equation}\label{n1}
\left[
\begin{array}{c}
\wt R(\lambda)  \\ \wt Q(\lambda)
\end{array}
\right]:=\Big(I_{2p}+\frac{i}{\lambda}C_{l+1}\Big)\times \ldots \times \Big(I_{2p}+\frac{i}{\lambda}C_{n}\Big) \left[
\begin{array}{c}
R(\lambda) \\ Q(\lambda)
\end{array}
\right]
\end{equation}
 satisfies condition 
\begin{equation}\label{n2}
W_{l+1}( \ov \lambda)^* \left[
\begin{array}{c}
\wt R(\lambda)  \\ \wt Q(\lambda)
\end{array}
\right]={\cal W}(\lambda) \left[
\begin{array}{c}
R(\lambda) \\ Q(\lambda)
\end{array}
\right].
\end{equation}
So if the pair $R$, $Q$ is admissible
for system (\ref{2.1}) on the interval $0 \leq k \leq n$, then the pair $\wt R$, $\wt Q$ is admissible
for system (\ref{2.1}) on the interval $0 \leq k \leq l$. Moreover Weyl function
$\vp$ of system (\ref{2.1}) on the interval $0 \leq k \leq n$ that is determined by $R$ and $Q$
coincides with the Weyl function
 of system (\ref{2.1}) on the interval $0 \leq k \leq l$ that is determined by $\wt R$ and $\wt Q$.
\end{Rk}
The next theorem solves inverse problem to recover system
(\ref{2.1}) from its Weyl  function.
\begin{Tm} \label{Tm2.2} Suppose $\vp$ is
a Weyl function of system (\ref{2.1})  satisfying conditions
(\ref{2.2}) and (\ref{2.3}).
 Then system (\ref{2.1}) is uniquely recovered
from the first $n+1$ Taylor coefficients $\{ \a_k \}_{k=0}^n$ of
$\, \displaystyle{\vp \left(i \Big( \frac{1+z}{1-z} \Big) \right)
}$ at $z=0$ by the following procedure.

First introduce $(n+1)p \times p$ matrices $\Phi_1$, $\Phi_2:$
\begin{equation}\label{2.11}
\Phi_1 = \left[
\begin{array}{c}
I_{p}  \\ I_p \\ \cdots \\ I_{p}
\end{array}
\right], \quad \Phi_2 = - \left[
\begin{array}{l}
\a_0  \\ \a_0+ \a_1 \\ \cdots \\ \a_0+ \a_1 + \ldots + \a_n
\end{array}
\right].
\end{equation}
Then introduce $(n+1)p \times 2p$ matrix $\Pi$ and $(n+1)p \times
(n+1)p$ block lower triangular matrix $A$ by their blocks $:$
$\quad \Pi=[\Phi_1 \quad \Phi_2]$,
\begin{equation}\label{2.12}
\begin{array}{lcllcl}
 A:=A(n) &=& \left\{ a_{j-k}^{\,} \right\}_{k,j=0}^n,
           & a_r  &=&  \left\{ \begin{array}{lll}
                                  0 \, & \mbox{ for }& r > 0   \\
                                 \displaystyle{\frac{i}{ {\, 2 \,}}} \,
                                 I_p
                                   & \mbox{ for }& r = 0   \\
                                 \, i \, I_p
                                   & \mbox{ for }& r < 0
                           \end{array} \right. \end{array}.
\end{equation}
Next we recover $(n+1)p \times (n+1)p$ matrix $S$ as a unique
solution of the matrix identity
\begin{equation} \label{2.13}
AS-SA^*=i \Pi \Pi^*.
\end{equation}
This solution is invertible and positive, i.e., $S>0$. Finally
matrices $\b(k)^*\b(k)$ are easily recovered from the formula
\begin{equation}\label{2.14}
\Pi^*S^{-1}\Pi=B^*B, \quad B: =B(n)= \left[
\begin{array}{c}
\b(0) \\ \b(1) \\ \cdots \\ \b(n)
\end{array}
\right].
\end{equation}
Now matrices $C_k$ and system (\ref{2.1}) are defined via
(\ref{2.2}).
\end{Tm}
\begin{proof}. Step 1. The method of the proof coincides with the
method of the proof of Theorem \ref{Tm1.2}. Put
\begin{equation}\label{2.15}
K(r) = \left[
\begin{array}{c}
K_0(r) \\ K_1(r) \\ \cdots \\ K_r(r)
\end{array}
\right],
\end{equation}
where $K_j(r)$ are $p \times (r+1)p$ matrices of the form
\begin{equation}\label{2.16}
K_j(r)=i \b(j)[\b(0)^* \ldots \b(j-1)^* \quad \b(j)^*/2 \quad 0
\ldots 0].
\end{equation}
From (\ref{2.14})-(\ref{2.16}) it follows that
\begin{equation}\label{2.17}
K(r)-K(r)^*=i B(r) B(r)^*.
\end{equation}
By induction we shall show in the next step that $K$ is similar to
$A$:
\begin{equation}\label{2.18}
K(r)=V_-(r)A(r)V_-(r)^{-1} \quad (0 \leq r \leq n),
\end{equation}
where $V_-(r)^{\pm 1}$ are block lower triangular matrices. Taking
into account (\ref{2.18}) and multiplying both sides  of
(\ref{2.17}) by $V_-(r)^{-1}$ from the left  and by $\big(V_-(r)^*
\big)^{-1}$ from the right  we get
\begin{equation}\label{2.19}
A(r)S(r)-S(r)A(r)^*=i \Pi(r) \Pi(r)^*,
\end{equation}
\begin{equation}\label{2.20}
S(r):=V_-(r)^{-1}\big(V_-(r)^* \big)^{-1}, \quad
\Pi(r):=V_-(r)^{-1}B(r).
\end{equation}
Moreover,  Step 3 will show that matrix $V_-(n)$ can be chosen so
that the equality
\begin{equation}\label{2.21}
\Pi=[\Phi_1 \quad \Phi_2]=V_-(n)^{-1}B(n)
\end{equation}
holds, i.e., $\Pi= \Pi(n)$. (Here $\Phi_1$ and $\Phi_2$ are given
by (\ref{2.11}).)

Identities (\ref{2.19}) have unique solutions $S(r)$ as the
spectra of $A$ and $A^*$ don't intersect. In particular, by
(\ref{2.13}) and (\ref{2.19}) one can see that $S:=S(n)$. Hence we
derive from (\ref{2.20}) and (\ref{2.21}) that $S>0$ and the first
equality in (\ref{2.14}) holds. It remains only to prove
(\ref{2.18}) and (\ref{2.21}).

Step 2. Now we shall consider block lower triangular matrices
$V_-(k)$ $(0 \leq k \leq n)$:
\begin{equation}\label{2.22}
V_-(0)=v_-(0)=\b_1(0), \quad V_-(k)= \left[
\begin{array}{cc}
V_-(k-1) & 0 \\ X(k) & v_-(k)
\end{array}
\right] \quad (k>0),
\end{equation}
where $v_-(l)$ are $p \times p$ matrices,  $X(k)=[X_0(k) \quad \wt
X(k) ]$ is $p \times k p$ matrix, $X_0(k)$ is an arbitrary $p
\times p$ block, and $\wt X(k)$, $v_-(k)$ are given by the
formulas
\[
\wt X(k)=i \Big( \b(k)[ \b(0)^* \ldots \b(k-1)^*]V_-(k-1)\left[
\begin{array}{c}
I_{(k-1)p}\\ 0
\end{array}
\right]-v_-(k)[I_p \ldots I_p] \Big)
\]
\begin{equation}\label{2.23}
\times \Big(A(k-2)+ \frac{i}{2} I_{(k-1)p}\Big)^{-1}, \quad
v_-(k)=\b(k)\b(k-1)^* v_-(k-1).
\end{equation}
According to (\ref{2.12}) we have $A(0)=(i/2)I_p$. From the second
relation in (\ref{2.2}) and definitions (\ref{2.15}) and
(\ref{2.16}) it is immediate that $K(0)=(i/2)I_p$ and so
(\ref{2.18}) is valid for $r=0$. Assume that (\ref{2.18}) is true
for $r=k-1$, and let us show that (\ref{2.18}) is true for $r=k$
too. It's easy to see that
\begin{equation}\label{2.25}
 V_-(k)^{-1}= \left[
\begin{array}{cc}
V_-(k-1)^{-1} & 0 \\ -v_-(k)^{-1} X(k)V_-(k-1)^{-1}  & v_-(k)^{-1}
\end{array}
\right].
\end{equation}
Then in view of definitions (\ref{2.12}) and (\ref{2.22}) our
assumption implies
\begin{equation}\label{2.26}
V_-(k)A(k)V_-(k)^{-1}=\left[
\begin{array}{cc}
K(k-1) & 0 \\ Y(k) & \frac{i}{2} I_{p}
\end{array}
\right],
\end{equation}
where $Y(k)=\Big[\big( X(k)A(k-1)+iv_-(k)[I_p \ldots I_p]\big)
\quad \, \frac{i}{2}v_-(k) \Big]$
\[
\times \left[
\begin{array}{c}
V_-(k-1)^{-1}  \\ -v_-(k)^{-1} X(k)V_-(k-1)^{-1}
\end{array}
\right].
\]
Rewrite product on the right hand side of the last formula as
\begin{equation}\label{2.27}
Y(k)=\Big( X(k)\big(A(k-1)- \frac{i}{2}I_{kp} \big)+iv_-(k)[I_p
\ldots I_p]\Big)V_-(k-1)^{-1}.
\end{equation}
From (\ref{2.12}) and (\ref{2.27}) it follows that
\begin{equation}\label{2.28}
Y(k)=\Big[ \Big(\wt X(k)\big(A(k-2)+ \frac{i}{2}I_{(k-1)p}
\big)+iv_-(k)[I_p \ldots I_p] \Big) \quad
iv_-(k)\Big]V_-(k-1)^{-1}.
\end{equation}
Notice that the sequence $[I_p \ldots I_p]$ of identity matrices
in (\ref{2.28}) is one block smaller than in (\ref{2.27}). By
(\ref{2.23}) and (\ref{2.28}) we have
\[
Y(k)=i \b(k) \Big[ [ \b(0)^* \ldots \b(k-1)^*]V_-(k-1)\left[
\begin{array}{c}
I_{(k-1)p}\\ 0
\end{array}
\right] \quad \, \b(k-1)^* v_-(k-1) \Big]
\]
\begin{equation}\label{2.29}
\times V_-(k-1)^{-1}.
\end{equation}
Finally formulas (\ref{2.22}) and (\ref{2.29}) imply
\begin{equation}\label{2.30}
Y(k)=i \b(k) [ \b(0)^* \ldots \b(k-1)^*].
\end{equation}
According to the second relation in (\ref{2.2}) and formulas
(\ref{2.16}) and (\ref{2.30}) we get
\begin{equation}\label{2.31}
\Big[Y(k) \quad  \frac{i}{2}I_{p} \Big]=K_k(k).
\end{equation}
Using now (\ref{2.15})  and (\ref{2.30}) one can see that the
right hand side of (\ref{2.26}) equals $K(k)$. Thus (\ref{2.18})
is true for $r=k$ and therefore for all $0 \leq r \leq n$.

Step 3. To derive (\ref{2.21}) we shall first prove that matrices
$V_-(r)$ given by (\ref{2.22}) and (\ref{2.23}) can be chosen so
that
\begin{equation}\label{2.32}
 V_-(r)^{-1}B_1(r)=\left[
\begin{array}{c}
I_{p}   \\ \cdots \\ I_{p}
\end{array}
\right], \quad B_1(r):= B(r)\left[
\begin{array}{c}
 I_p \\0
\end{array}
\right]=\left[
\begin{array}{c}
 \b_1(0) \\ \cdots \\   \b_1(r)
\end{array}
\right].
\end{equation}
In other words the arbitrary till now  blocks $X_0(r)$ can be
chosen so. Indeed, by definition in (\ref{2.14}) and the first
equality in (\ref{2.22}) formula (\ref{2.32}) is true for $r=0$.
Assume that (\ref{2.32}) is true for $r=k-1$. Then from
(\ref{2.25}) it follows that (\ref{2.32}) is true for $r=k$ if
\begin{equation}\label{2.33}
-v_-(k)^{-1} X(k)\left[
\begin{array}{c}
I_{p}   \\ \cdots \\ I_{p}
\end{array}
\right]+v_-(k)^{-1}\b_1(k)=I_p.
\end{equation}
It implies that we get equality (\ref{2.32})  for $r=k$ putting
\begin{equation}\label{2.34}
X_0(k)=\b_1(k)-v_-(k)-\wt X(k)\left[
\begin{array}{c}
I_{p}   \\ \cdots \\ I_{p}
\end{array}
\right].
\end{equation}
Hence by a proper choice of matrices $X_0(r)$ we  obtain
(\ref{2.32}) for all $r \leq n$.

It remains to prove that
\begin{equation}\label{2.35}
 V_-(n)^{-1}B_2(n)=\Phi_2, \quad B_2(n):=\left[
\begin{array}{c}
 \b_2(0) \\ \cdots \\   \b_2(n)
\end{array}
\right].
\end{equation}
For that purpose we shall consider matrix function $W_{n+1}(\lambda)$,
that is used in (\ref{2.4}) to define coefficients of the M\"obius
transformation (\ref{2.5}) (and thus Weyl functions). Namely we
shall prove so called transfer matrix function representation of
$W_{n+1}(\lambda)$:
\begin{equation}\label{2.36}
W_{n+1}(\lambda)=\left(\frac{\lambda-i}{\lambda}\right)^{n+1}w_A \Big(n,
\frac{\lambda}{2}\Big),
\end{equation}
where
\begin{equation}\label{2.37}
w_A(r, \lambda)=I_{2p}-i \Pi(r)^*S(r)^{-1}\big(A(r)- \lambda I_{(r+1)p}
\big)^{-1} \Pi(r).
\end{equation}
Transfer matrix functions of the form (\ref{2.37}) have been
introduced and studied by Lev Sakhnovich \cite{SaL0}-\cite{SaL3}.
In particular, identity (\ref{2.19}) implies
\[
w_A(r, \mu)^*w_A(r,\lambda)=I_{2p}
\]
\begin{equation}\label{2.38}
+i(\ov \mu - \lambda) \Pi(r)^*\big(A(r)^*-\ov \mu
I_{(r+1)p}\big)^{-1}S(r)^{-1}\big(A(r) - \lambda
I_{(r+1)p}\big)^{-1}\Pi(r).
\end{equation}
Moreover, according to factorization Theorem 4 from \cite{SaL0}
(see also \cite{SaL2}, p. 188) we have
\[
w_A(r, \lambda)=\Big(I_{2p} -i \Pi(r)^*S(r)^{-1}P^*\big(PA(r)P^*- \lambda
I_p \big)^{-1}\big(PS(r)^{-1}P^*\big)^{-1}
\]
\begin{equation}\label{2.39}
\times P S(r)^{-1}\Pi(r) \Big) w_A(r-1, \lambda), \quad P=[0 \quad
\ldots \quad 0 \quad I_p].
\end{equation}
Taking into account (\ref{2.12}), (\ref{2.20}), and (\ref{2.22})
we obtain
\begin{equation}\label{2.40}
\big(PA(r)P^*- \lambda I_p \big)^{-1}=(\frac{i}{2}- \lambda )^{-1}I_p, \quad
PS(r)^{-1}P^*=v_-(r)^*v_-(r),
\end{equation}
\begin{equation}\label{2.41}
P S(r)^{-1}\Pi(r)=v_-(r)^*P B(r)=v_-(r)^* \b(r).
\end{equation}
Substitute (\ref{2.40}) and (\ref{2.41}) into (\ref{2.39}) to get
\begin{equation}\label{2.42}
w_A(r, \frac{\lambda}{2})=\Big(I_{2p} -\frac{2i}{i- \lambda}\b(r)^*\b(r)
\Big)w_A(r-1, \frac{\lambda}{2}).
\end{equation}
From the definitions (\ref{2.12}), (\ref{2.20}), and (\ref{2.37})
we also easily derive
\begin{equation}\label{2.43}
w_A(0, \frac{\lambda}{2})=I_{2p} -\frac{2i}{i- \lambda}B(0)^*B(0)=I_{2p}
-\frac{2i}{i- \lambda}\b(0)^*\b(0).
\end{equation}
On the other hand  system (\ref{2.1}) with additional condition
(\ref{2.2}) can be rewritten as
\begin{equation}\label{2.44}
W(r+1, \lambda)=\frac{\lambda - i}{\lambda}\Big(I_{2p}-\frac{2i}{i-
\lambda}\b(r)^*\b(r) \Big)W(r, \lambda).
\end{equation}
In view of the normalization $W(0)=I_{2p}$ formulas
(\ref{2.42})-(\ref{2.44}) imply  (\ref{2.36}).

From (\ref{2.36}) and (\ref{2.38}) it follows that
\begin{equation}\label{2.45}
W(n+1, \lambda)W(n+1, \ov
\lambda)^*=\left(\frac{\lambda-i}{\lambda}\right)^{n+1}\left(\frac{\lambda+i}{\lambda}\right)^{n+1}.
\end{equation}
Let us include now Weyl functions into consideration and put
\begin{equation}\label{2.46}
{\cal A}(\lambda):=\left| \frac{\lambda}{\lambda-i} \right|^{2n+2}[\vp(\lambda)^*
\quad I_p]W(n+1, \lambda)^*W(n+1, \lambda)\left[
\begin{array}{c}
\vp(\lambda) \\ I_{p}
\end{array}
\right].
\end{equation}
According to (\ref{2.5}), (\ref{2.4}), and (\ref{2.45}) we have
\[
{\cal A}(\lambda)=\left| \frac{\lambda+i}{\lambda} \right|^{2n+2} \bigl(
\big({\cal W}_{21}(\lambda ) R(\lambda )+{\cal W}_{22}(\lambda )Q(\lambda
)\big)^*\bigr)^{-1}
\]
\begin{equation}\label{2.47}
\times \Big(R(\lambda)^*R(\lambda)+ Q(\lambda)^*Q(\lambda)\Big)\bigl( {\cal W}_{21}(\lambda
) R(\lambda )+{\cal W}_{22}(\lambda )Q(\lambda )\bigr)^{-1}.
\end{equation}
By  (\ref{2.6}) and (\ref{2.47})  ${\cal A}$ is bounded in the
neighborhood of $\lambda=i$:
\begin{equation}\label{2.48}
\|{\cal A}(\lambda)\|=O(1) \quad {\mathrm{for}} \quad \lambda \to i.
\end{equation}
Substitute now (\ref{2.36}) and (\ref{2.38}) into (\ref{2.46}) to
obtain
\[
{\cal A}(\lambda)=[\vp(\lambda)^* \quad I_p] \Big( I_{2p}+\frac{i}{2}(\ov \lambda
- \lambda) \Pi(n)^*\big(A(n)^*-\frac{\ov \lambda}{2}
I_{(n+1)p}\big)^{-1}S(n)^{-1}
\]
\begin{equation}\label{2.49}
\times \big(A(n) - \frac{\lambda}{2}
I_{(n+1)p}\big)^{-1}\Pi(n)\Big)\left[
\begin{array}{c}
\vp(\lambda) \\ I_{p}
\end{array}
\right].
\end{equation}
Notice that $S(n)>0$. Hence formulas (\ref{2.48}) and (\ref{2.49})
imply that
\begin{equation}\label{2.50}
\left\| \big(A(n) - \frac{\lambda}{2} I_{(n+1)p}\big)^{-1}\Pi(n)\left[
\begin{array}{c}
\vp(\lambda) \\ I_{p}
\end{array}
\right]\right\|=O(1) \quad {\mathrm{for}} \quad \lambda \to i.
\end{equation}
Recall that $\Pi(n)=V_-(n)^{-1}B(n)$ and $A(n)$ is denoted by $A$.
Represent now $\Pi(n)$ in the block form
\begin{equation}\label{2.50'}
\Pi(n)=[\Phi_1(n) \quad \Phi_2(n)], \quad
\Phi_k(n)=V_-(n)^{-1}B_k(n) \quad (k=1,2).
\end{equation}
According to (\ref{2.11}) and (\ref{2.32}) we have
$\Phi_1(n)=\Phi_1$. Hence multiplying the matrix function on the
left hand side of (\ref{2.50}) by $\Big(\Phi_1^*\big(A -
\frac{\lambda}{2} I_{(n+1)p}\big)^{-1}\Phi_1 \Big)^{-1}\Phi_1^*$ we
derive
\[\left\|
\vp(\lambda)+\Big(\Phi_1^*\big(A - \frac{\lambda}{2}
I_{(n+1)p}\big)^{-1}\Phi_1 \Big)^{-1}\Phi_1^*\big(A - \frac{\lambda}{2}
I_{(n+1)p}\big)^{-1}\Phi_2(n)\right\|
\]
\begin{equation}\label{2.51}
=O\left(\left\|\Big(\Phi_1^*\big(A - \frac{\lambda}{2}
I_{(n+1)p}\big)^{-1}\Phi_1 \Big)^{-1}\right\|\right) \quad
{\mathrm{for}} \quad \lambda \to i.
\end{equation}
The matrix $A - \frac{\lambda}{2} I_{(n+1)p}$ is easily inverted
explicitly (see, for instance, formula (1.10) in \cite{SaAtepl}).
As a result one obtains
\begin{equation}\label{2.52}
\Phi_1^*\big( A - \frac{\lambda}{2}
I_{(n+1)p}\big)^{-1}=\frac{2}{i-\lambda}[q^n \quad q^{n-1} \quad \ldots
\quad q \quad I_p ], \quad q:=\frac{\lambda +i}{\lambda -i}.
\end{equation}
Moreover we get
\begin{equation}\label{2.53}
\Phi_1^*\big( A - \frac{\lambda}{2} I_{(n+1)p}\big)^{-1}\Phi_1=\frac{2}{i-\lambda}\Big(q^{n+1}-I_p \Big)\Big(q-I_p \Big)^{-1}
\end{equation}
 Putting $\displaystyle{\lambda =i \Big( \frac{1+z}{1-z}
\Big)}$, i.e., $\displaystyle{z =\Big( \frac{\lambda -i}{\lambda+i}
\Big)}$ we derive from (\ref{2.53}) that
\begin{equation}\label{2.53'}
\Big(\Phi_1^*\big(A - \frac{\lambda}{2}
I_{(n+1)p}\big)^{-1}\Phi_1 \Big)^{-1}=\big(-iz^{n+1}+O(z^{2n+2})\big)I_p \quad
(z \to 0).
\end{equation}
Taking into account (\ref{2.52}) and (\ref{2.53'}) we
rewrite (\ref{2.51}) as 
\begin{equation}\label{2.54}
\left\| \vp\Big(i \Big( \frac{1+z}{1-z}
\Big)\Big)+(1-z)[I_p \quad zI_p \quad z^2I_p \quad \ldots]\Phi_2(n)\right\|
=O(z^{n+1})  
\end{equation}
for $z \to 0$. From (\ref{2.11}) and (\ref{2.54}) follows that $\Phi_2(n)=\Phi_2$, i.e., (\ref{2.35})
is true.
As $\Phi_1(n)=\Phi_1$ and $\Phi_2(n)=\Phi_2$, so $\Pi(n)=\Pi$ and formula (\ref{2.21})
is finally proved. 
\end{proof}
From Theorem \ref{Tm2.2} and Remark  \ref{Rk2.1'}  follows  Borg-Marchenko type result.
\begin{Tm}\label{Tm2.4}
Let $\wt \varphi$ and $\wh \varphi$ be  Weyl functions of the two
discrete Dirac type systems (\ref{2.1}) that satisfy conditions
(\ref{2.2}) and (\ref{2.3}). Denote by $\wt C_k$ ($0 \leq k \leq \wt n$) the potentials $C_k$ of  the first system and by  
$\wh C_k$ ($0 \leq k \leq \wh n$) the potentials
of the second system. Denote Taylor coefficients of $\wt \vp\Big(i \Big( \frac{1+z}{1-z}
\Big)\Big)$ and $\wh \vp\Big(i \Big( \frac{1+z}{1-z}
\Big)\Big)$ at $z=0$ by $\{\wt \a_k \}$ and $\{\wh \a_k \}$, respectively, and assume that
$\wt \a_k=\wh \a_k$ for $k \leq l$ ($l \leq \wt n$, $l \leq \wh n$). Then we have
$\wt C_k=\wh C_k$ for $k \leq l$.
\end{Tm}
\begin{proof}.  According to Remark  \ref{Rk2.1'} $\wt \varphi$ and $\wh \varphi$ are Weyl functions
of the first and second systems, respectively, on the interval $0 \leq k \leq l$. By Theorem \ref{Tm2.2}
these systems on the interval $0 \leq k \leq l$ are uniquely recovered by the first $l+1$
Taylor coefficients of the Weyl functions.
\end{proof}
Step 3 of the proof of Theorem  \ref{Tm2.2} implies the following corollary.
\begin{Cy}\label{Cy2.3}
Weyl functions of system (\ref{2.1}) that satisfies conditions
(\ref{2.2}) and (\ref{2.3}) admit Taylor representation
\begin{equation}\label{2.55}
 \vp\Big(i \Big( \frac{1+z}{1-z}
\Big)\Big)=-\psi_0+(\psi_0-\psi_1)z+\ldots+(\psi_{n-1}-\psi_{n})z^n+O(z^{n+1}) \quad (z \to 0),
\end{equation}
where $\psi_k$ are $p \times p$ blocks of $\Phi_2(n)=\{ \psi_k \}_{k=0}^n$. Here $\Phi_2(n)$
is given by the second equality in (\ref{2.50'}), where $V_-(n)$ is defined by formulas (\ref{2.22}),
(\ref{2.23}), and (\ref{2.34}).
\end{Cy}
Moreover, from the proof of Theorem  \ref{Tm2.2} follows a complete description of the Weyl functions
in terms of the Taylor coefficients.
\begin{Tm}\label{Tm2.5}
(i) Let system (\ref{2.1}) be given on the interval $0 \leq k \leq n$
and satisfy (\ref{2.2}), (\ref{2.3}). Then analytic  at $\lambda =i$ matrix function  $\vp$
is a  Weyl function  of this system if and only if it admits
expansion  (\ref{2.55}), where matrices $\psi_k$ are defined in Corollary \ref{Cy2.3}. \\
(ii) Suppose $\vp$ is a $p \times p$ matrix function analytic at $\lambda =i$. Then $\vp$
is a Weyl function of some system (\ref{2.1}) given on the interval $0 \leq k \leq n$
and satisfying (\ref{2.2}), (\ref{2.3}) if and only if the matrix $S$ uniquely defined by the
identity $AS-SA^*=i \Pi \Pi^*$ is invertible. Here $\Pi=[\Phi_1 \quad \Phi_2]$ is given by
(\ref{2.11}), where $\{\a_k\}$ are Taylor coefficients of $ \vp\Big(i \Big( \frac{1+z}{1-z}
\Big)\Big)$ at $z=0$.
\end{Tm}
\begin{proof}.
Suppose $\vp$ is a $p \times p$ matrix function analytical at $\lambda =i$ and satisfying  (\ref{2.55}),
and put
\begin{equation}   \label{2.55'}
 \left[
\begin{array}{c}
R(\lambda) \\ Q(\lambda)
\end{array}
\right]:= { \cal W }(\lambda)^{-1}  \left[
\begin{array}{c}
\vp(\lambda) \\ I_p
\end{array}
\right].
\end{equation}
Fix a Weyl function $\wt \vp$ of   system (\ref{2.1}), and denote by $\wt R$, $\wt Q$
some admissible pair that grants representation (\ref{2.5}) of $\wt \vp$.
Rewrite (\ref{2.55'}) in the form
\begin{equation}   \label{2.56}
 \left[
\begin{array}{c}
R(\lambda) \\ Q(\lambda)
\end{array}
\right]={\cal W}(\lambda)^{-1}   \left[
\begin{array}{c}
\wt \vp(\lambda) \\ I_p
\end{array}
\right]+ {\cal W}(\lambda)^{-1} \left[
\begin{array}{c}
\vp(\lambda)- \wt \vp(\lambda) \\ 0
\end{array}
\right].
\end{equation}
As $\wt \vp$ is the M\"obius transformation (\ref{2.5}) of
the admissible pair $\wt R$, $\wt Q$, we have
\begin{equation}   \label{2.57}
{\cal W}(\lambda)^{-1}   \left[
\begin{array}{c}
\wt \vp(\lambda) \\ I_p
\end{array}
\right]= \left[
\begin{array}{c}
\wt R(\lambda) \\ \wt Q(\lambda)
\end{array}
\right]\bigl( {\cal W}_{21}(\lambda )\wt R(\lambda )+{\cal W}_{22}(\lambda
) \wt Q(\lambda )\bigr)^{-1}.
\end{equation}
Thus the first summand on the right-hand side of (\ref{2.56}) is
analytic at $\lambda=i$. Taking into account that expansion
(\ref{2.55}) is valid for $\vp$ and $\wt \vp$ we derive
$\vp(\lambda)- \wt \vp(\lambda)=O\big((\lambda-i)^{n+1}\big)$ for
$\lambda \to i$. From (\ref{2.45}) it follows also that
$\displaystyle {\cal
W}(\lambda)^{-1}=\frac{\lambda^{2n+2}}{(\lambda+i)^{n+1}(\lambda-i)^{n+1}}W(n+1,
\lambda)$. So the second summand on the right-hand side of
(\ref{2.56}) is analytic at $\lambda=i$ too. Therefore the pair
$R$, $Q$ is analytic at $\lambda=i$.  Moreover, according to
(\ref{2.55'}) we have ${\cal W}_{21}(\lambda ) R(\lambda )+{\cal
W}_{22}(\lambda)Q(\lambda)=I_p$. Hence, inequality (\ref{2.6})
holds and the pair $R$, $Q$ is admissible. It easily follows from
(\ref{2.55'}) that $\vp$ admits representation (\ref{2.5}) with
this $R$, $Q$, i.e., $\vp$ is a Weyl function of our system.

Vice versa, if $\vp$ is a Weyl function of our system, then by Corollary \ref{Cy2.3} 
the expansion (\ref{2.55}) is true. The statement (i) is proved.

According to the proof  of Theorem \ref{Tm2.2} when $\vp$ is a Weyl function, then matrix $S$
uniquely defined by the identity $AS-SA^*=i \Pi \Pi^*$ is invertible. It remains to show, that
if $S$ is invertible, then $\vp$ is a Weyl function. For this purpose introduce notations
\begin{equation}   \label{2.58}
P_r=[0 \quad
\ldots \quad 0 \quad I_p], \quad \wt P_r=[I_{(r+1)p}  \quad 0 ], \quad S(r)= \wt P_r S  \wt P_r^*,
\quad \Pi(r)= \wt P_r \Pi,
\end{equation}
where $P_r$ are $p \times (r+1)p$ and $\wt P_r$ are  $(r+1)p \times (n+1)p$ matrices. Assume
that $\det \, S \not=0$. As the identity (\ref{2.13}) is equivalent to 
\[
S(A^*-\lambda I_{(n+1)p})^{-1}-(A-\lambda I_{(n+1)p})^{-1}S=i(A-\lambda I_{(n+1)p})^{-1}\Pi \Pi^* (A^*-\lambda I_{(n+1)p})^{-1},
\]
using residues we standardly derive
\begin{equation}   \label{2.59}
S=\frac{1}{2 \pi} \int_{- \infty}^{ \infty}(A-\lambda I_{(n+1)p})^{-1}\Pi \Pi^* (A^*-\lambda I_{(n+1)p})^{-1}d \lambda.
\end{equation}
Therefore from $\det \, S \not=0$ it follows that $S>0$, and so we have $S(r)>0$ and $\det \, S(r) \not=0$. Put
\begin{equation}   \label{2.60}
\b(r):=\big(P_rS(r)^{-1}P_r \big)^{-\frac{1}{2}}P_rS(r)^{-1}\Pi(r) \quad (0 \leq r \leq n).
\end{equation}
Matrices $\b(r)$ satisfy conditions (\ref{2.2}) and (\ref{2.3}). Indeed, from  (\ref{2.60})  we get
\begin{equation}   \label{2.61}
\b(r)\b(r)^*=\big(P_rS(r)^{-1}P_r \big)^{-\frac{1}{2}}P_rS(r)^{-1}\Pi(r)\Pi(r)^*
S(r)^{-1}P_r^*\big(P_rS(r)^{-1}P_r \big)^{-\frac{1}{2}}.
\end{equation}
From (\ref{2.13}) follows identity
(\ref{2.19}), i.e., 
\[
S(r)^{-1} \Pi(r)\Pi(r)^* S(r)^{-1}=-i(S(r)^{-1}A(r)-A(r)^*S(r)^{-1}). 
\]
It is also true that $P_r A(r)^*=- \frac{i}{2}P_r$ and $A(r)P_r^*= \frac{i}{2}P_r$. Hence we can rewrite (\ref{2.61})
as
\[
\b(r)\b(r)^*=-i\big(P_rS(r)^{-1}P_r \big)^{-\frac{1}{2}}P_r\Big(S(r)^{-1}A(r)
-A(r)^*S(r)^{-1}\Big)
\]
\begin{equation}   \label{2.62}
\times P_r^*\big(P_rS(r)^{-1}P_r \big)^{-\frac{1}{2}}=I_p.
\end{equation}
Thus the condition from (\ref{2.2}) on $\b$  is proved. Notice further that according to (\ref{2.11}) and 
(\ref{2.60}) we have
\begin{equation}   \label{2.63}
\b_1(0)=\big(S(0)^{-1} \big)^{-\frac{1}{2}}S(0)^{-1}=S(0)^{-\frac{1}{2}},
\end{equation}
and so the first inequality in (\ref{2.3}) is also true. To prove  the second inequality in (\ref{2.3})
we should write down  $S(r)^{-1}$ $(r>0)$ in the block form:
\begin{equation}   \label{2.64}
 S(r)^{-1}=\left[
\begin{array}{lr}
S(r-1)^{-1}+S(r-1)^{-1}S_{12}tS_{21}S(r-1)^{-1} & -S(r-1)^{-1}S_{12}t\\ -tS_{21}S(r-1)^{-1} & t
\end{array}
\right],
\end{equation}
where $t=\Big(s-S_{21}S(r-1)^{-1}S_{12}\Big)^{-1}$, and $S_{12}$, $S_{21}$ and $s$ are blocks of $S(r)$:
\begin{equation}   \label{2.65}
S(r)=\left[
\begin{array}{lr}
S(r-1) & S_{12} \\ S_{21} & s
\end{array}
\right].
\end{equation}
Formula (\ref{2.64}) is easily checked directly. (Conditions  $\det \, S(r) \not=0$ and $\det \, S(r-1) \not=0$
imply invertibility of $s-S_{21}S(r-1)^{-1}S_{12}$; see, for instance, \cite{SaL3}, p. 21.) In view of (\ref{2.60})
and (\ref{2.64})  the second inequality in (\ref{2.3}) is equivalent to
\begin{equation}   \label{2.66}
\det \Big( [-S_{21}S(r-1)^{-1} \quad I_p]\Pi(r)\Pi(r-1)^*S(r-1)^{-1}P_{r-1}^*\Big) \not=0 \quad (r>0).
\end{equation}
Recall that $S(r-1)$ also satisfies operator identity: 
\[
A(r-1)S(r-1)-S(r-1)A(r-1)^*=i\Pi(r-1)\Pi(r-1)^*. 
\]
Therefore similar to the proof of (\ref{2.62}) one gets
\[
 [-S_{21}S(r-1)^{-1} \quad I_p]\Pi(r)\Pi(r-1)^*S(r-1)^{-1}P_{r-1}^* =-\frac{1}{2}S_{21}S(r-1)^{-1}P_{r-1}^* 
\]
\begin{equation}   \label{2.67}
 -iS_{21}A(r-1)^*S(r-1)^{-1}P_{r-1}^*+[I_p \quad \psi_r]\Pi(r-1)^*S(r-1)^{-1}P_{r-1}^*,
\end{equation}
where $\{\psi_k\}$ are the blocks of $\Phi_2$. Consider now the first $r$ blocks in the lower
block rows on the both sides of  the identity (\ref{2.19}). In view of (\ref{2.65}) we have
\[
\frac{i}{2}S_{21}+i[I_p \quad \ldots \quad I_p]S(r-1)-S_{21}A(r-1)^*=i[I_p \quad \psi_r]\Pi(r-1)^*,
\quad {\mathrm{i.e.,}}
\]
\begin{equation}   \label{2.68}
i\Big(\frac{i}{2}S_{21}-S_{21}A(r-1)^*-i[I_p \quad \psi_r]\Pi(r-1)^*\Big)S(r-1)^{-1}P_{r-1}^*=I_p.
\end{equation}
By (\ref{2.67}) and (\ref{2.68})  inequality (\ref{2.66})  holds. 

So formula (\ref{2.60}) defines
system (\ref{2.1}) that satisfies conditions (\ref{2.2}) and (\ref{2.3}).
Moreover, by Theorem 4 \cite{SaL0} the transfer matrix functions $w_A(r,\lambda)$ corresponding
to our matrices $S(r)$ admit factorizations (\ref{2.39}). In view of (\ref{2.60}) this implies
\begin{equation}   \label{2.69}
w_A(n, \frac{\lambda}{2})=\prod_{r=0}^n\Big(I_{2p} -\frac{2i}{i- \lambda}\b(r)^*\b(r)
\Big).
\end{equation}
Compare (\ref{2.69}) with (\ref{2.44}) to derive for the fundamental solution of the constructed system
the equality $W(n+1, \lambda)=\Big(\frac{\lambda -i}{\lambda}\Big)^{n+1}w_A(n, \frac{\lambda}{2})$, i.e.,
\begin{equation}   \label{2.70}
{\cal W}(\lambda)^{-1}=\Big(\frac{\lambda }{\lambda+i}\Big)^{n+1}w_A(n, \frac{\lambda}{2}).
\end{equation}

Let us consider now $\vp(\lambda)$ and prove (\ref{2.50}). It is easy to see (\cite{SaAtepl}, p. 452) that
the $k$-th block row of  $\big(A- \frac{\lambda}{2} I_{(n+1)p}\big)^{-1}$ is given by the formula
\[
T_k(\lambda)=\frac{i(1-z)}{z^{k+1}}
\]
\begin{equation}   \label{2.71}
\times \Big[(1-z)I_p \quad (1-z)zI_p \quad \ldots \quad (1-z)z^{k-1}I_p \quad z^kI_p \quad 0 \quad \ldots
\quad 0 \Big],
\end{equation}
where $\displaystyle{z =\Big( \frac{\lambda -i}{\lambda+i}\Big)}$. It follows also that 
\begin{equation}   \label{2.72}
\big(A-  \frac{\lambda}{2} I_{(n+1)p}\big)^{-1}\Phi_1= \frac{i(1-z)}{z}\,  {\mathrm  {col}}\Big[I_p \quad z^{-1}I_p \quad
\ldots \quad z^{-n}I_p \Big],
\end{equation}
where col means column. By (\ref{2.71}) and (\ref{2.72}) we have
\[
\big(A - \frac{\lambda}{2} I_{(n+1)p}\big)^{-1}\Pi \left[
\begin{array}{c}
\vp(\lambda) \\ I_{p}
\end{array}
\right]
\]
\begin{equation}   \label{2.73}
=\left\{ \frac{i(1-z)}{z^{k+1}}\Big( \vp\Big(i \Big( \frac{1+z}{1-z}
\Big)\Big)+\psi_0+(\psi_1-\psi_0)z+\ldots+(\psi_{k}-\psi_{k-1})z^k \Big) \right\}_{k=0}^n.
\end{equation}
Formulas (\ref{2.55}) and (\ref{2.73}) imply (\ref{2.50}). From (\ref{2.37}),  (\ref{2.50}),
and (\ref{2.70}) it follows that $\displaystyle \left\| { \cal W }(\lambda)^{-1}  \left[
\begin{array}{c} \vp(\lambda) \\ I_p \end{array} \right]\right\|=O(1) $ for $\lambda \to i$. Thus
the pair $R$, $Q$ given by (\ref{2.55'}) is admissible. As $\vp$ admits representation
(\ref{2.5}) with $R$, $Q$ given by (\ref{2.55'}), so $\vp$ is a Weyl function of the
constructed system.
\end{proof}
\begin{Ee} \label{Een2}
According to (\ref{e5}) for system considered in Example \ref{Een1} we have $\a_0=0$ and $\a_1=I_p$.
Hence, by Theorem \ref{Tm2.5} the set of Weyl functions $\vp$ of this system is defined by the expansion
\begin{equation}\label{e6}
 \vp\Big(i \Big( \frac{1+z}{1-z}
\Big)\Big)=z+O(z^2), \quad z \to 0.
\end{equation}
\end{Ee}

\end{document}